\documentclass[12pt,leqno]{amsart}
\usepackage{amsmath}
\usepackage{graphicx,color}
\newtheorem{theorem}{Theorem}[section]

\newtheorem{lemma}[theorem]{Lemma}
\newtheorem{proposition}[theorem]{Proposition}
\theoremstyle{definition}
\newtheorem{definition}[theorem]{Definition}
\theoremstyle{remark}

\numberwithin{equation}{section}

\def\R{\mathbb{R}} %Real R
\def\C{\mathbb{C}} %Complex C
\def\D{\mathbb {D}}
\def\Re{{ Re\,}}
\def\Im{{ Im\,}}

\begin{document}

\sloppy
%10000

\title{Levi-flat filling of real two-spheres in symplectic manifolds (I)}

\author{Herv{\'e} Gaussier{*} and Alexandre Sukhov{**}}

%\address{}


\date{\number\year-\number\month-\number\day}

\maketitle

{\small
* Universit\'e Joseph Fourier, 100 rue des Maths, 38402 Saint Martin d'H\`eres, France, herve.gaussier@ujf-grenoble.fr

** Universit\'e des Sciences et Technologies de Lille, Laboratoire
Paul Painlev\'e,
U.F.R. de
Math\'e-matique, 59655 Villeneuve d'Ascq, Cedex, France,
 sukhov@math.univ-lille1.fr
}

\bigskip
Abstract. Let $(M,J,\omega)$ be a manifold with an almost complex structure $J$ tamed by a symplectic form $\omega$. We suppose that $M$ has the complex dimension two, is Levi convex and with bounded geometry. We prove that a real two-sphere with two elliptic points, embedded into the boundary of $M$ can be foliated by the boundaries of pseudoholomorphic discs.

\bigskip

R\'esum\'e. Soit $(M,J,\omega)$ une vari\'et\'e dont la structure presque complexe $J$ est tamed par la forme symplectique $\omega$. On suppose $M$ de dimension complexe deux, Levi convexe et \`a g\'eom\'etrie born\'ee. On d\'emontre que toute 2-sph\`ere poss\'edant deux points elliptiques et plong\'ee dans le bord de $M$ est feuillet\'ee par des bords de disques pseudoholomorphes.

\bigskip

\section{Introduction}
This expository paper concerns the problem of filling real two-spheres in almost complex manifolds by pseudoholomorphic discs. We prove the following statement :
\begin{theorem}
\label{elliptic-thm}
Let $(M,J,\omega)$ be an almost complex manifold of complex dimension 2 with a
taming symplectic form, with bounded geometry. Assume that $M$ does not contain any compact $J$-complex
sphere and that the boundary $\partial M$ of $M$ is a smooth
Levi convex hypersurface that does not contain any germ of a nonconstant $J$-holomorphic disc. Let $S^2$ be a real $2$-sphere with two
elliptic points, embedded into $\partial M$. Then there exists a unique smooth one parameter family of
disjoint Bishop discs for $S^2$ filling a real Levi-flat hypersurface $\Gamma
\subset M$ with boundary $S^2$. This hypersurface is smooth up to the boundary
except at the two elliptic points.
\end{theorem}

Though results of this type admit important applications and have been obtained by several authors in various forms, it seems difficult to find a reference for a complete 
proof of the above statement. This is our first motivation. The above result and the techniques of the present work will be used in a forthcoming paper
devoted the the problem of filling two-spheres with elliptic and hyperbolic points. The present work allows us to focus the second paper on properties of discs near hyperbolic points.

Recall that according to classical results a real two-sphere does not admit a totally real embedding into a two dimensional complex manifold : any embedding of $S^2$ has complex points. In the generic case these points are isolated and are either elliptic or hyperbolic. Parabolic points can be removed by a  small deformation. The first related local result goes back to E.Bishop \cite{bi} who proved the existence of a family of holomorphic discs with boundaries attached to a neigborhood of an analytic point in a real two dimensional sphere in $\C^2$.
The global result is due to E.Bedford-B.Gaveau \cite{be-ga}. They proved the following. Let $\Omega$ be a strictly pseudoconvex domain in $\C^2$ and $S^2$ be a real 2-sphere embedded into $\partial \Omega$. If $S^2$ has exactly two elliptic points $p,\ q$ and is totally real outside these points then $S^2$ is the boundary of a Levi-flat hypersurface $\Gamma$ foliated by holomorphic discs with boundaries on $S^2 \backslash \{p,q\}$.
The hypersurface $\Gamma$ can be viewed as a resolution of the Plateau problem and has other important properties related to polynomially convex hulls and envelopes of holomorphy. This result was extended by M.Gromov \cite{gr} to the almost complex case under the assumption that the almost complex structure is real analytic and integrable near $S^2$. This extra assumption was removed by R.Ye in the paper \cite{ye}.
These results have many important consequences in complex and contact geometry. One of them is the fundamental observation due to M.Gromov that the characteristic foliation on any 2-sphere embedded into the 3-sphere with the standard contact structure has no closed leaf. Another important application is due to H.Hofer \cite{ho} for his proof of the Weinstein's conjecture for contact 3-manifolds $M$ with $\pi_2(M) \neq \{ 0 \}$. Y.Eliashberg discovered several applications in symplectic topology \cite{el}.

The above mentioned results were obtained for 2-spheres contained in the boundary of a strictly pseudoconvex domain. In the present paper we consider the filling of a 2-sphere embedded to the boundary of a pseudoconvex domain assuming that this boundary contains no non-constant holomorphic discs. This condition is automatically satisfied for strictly pseudoconvex domains so the class of domains under consideration is substantially larger. Y.Eliashberg-W.Thurston \cite{el-th} introduced a new class of structures called confoliations, intermediate objects between contact structures and foliations. One of the most important examples of such structures is provided by the holomorphic tangent bundle of a weakly pseudoconvex hypersurface in $\C^2$. This can be viewed as a contact structure with degeneracies. The problem of filling a real 2-sphere in a confoliated 3-manifold by  pseudo-holomorphic discs   was studied by R.Hind \cite{hi}. His approach is based on a result due to E.Bedford-W.Klingenberg \cite{be-kl} on the Levi-flat filling of spheres in the presence of hyperbolic points and on a theorem of Y.Eliashberg-W.Thurston
on the approximation of a confoliation by contact structures. However several
important  arguments of the proof are dropped in \cite{hi}. 
The present paper  presents a detailed proof in the elliptic case, based on a
quite different approach. Our main idea is to use an almost complex version of
the well-known result of Diederich-Fornaess \cite{di-fo} on the existence of
bounded exhaustion plurisubharmonic functions for weakly pseudoconvex domains
in $\C^n$. The corresponding almost complex version is obtained in
\cite{di-su}. It removes many technical problems in the proof. In particular
this approach allows to give a direct proof without using
Eliashberg-Thurston's theory. The results of \cite{di-su} are valid in any
dimension in contrast to the Eliashberg-Thurston theory only working in complex
dimension two. This allows to apply our techniques in the study of
confoliations in the higher dimensional case for some classes of complex structures. At the end of the paper (see Remark 3 in the last section) we indicate how to use our approach in order to obtain new results on the filling of two-spheres in weakly pseudoconvex boundaries by holomorphic discs.
We also point out that the examples due to Y.Eliashberg \cite{el} and J.E.Fornaess-D.Ma \cite{fo-ma} show that the condition of pseudoconvexity of $\partial M$ in Theorem~\ref{elliptic-thm} 
cannot be dropped.

Keeping in mind forthcoming papers and for convenience of the reader we give a reasonably self contained exposition.

%The aim of this paper is to prove the following Theorem :

%\begin{theorem}
%\label{main-thm}
%Let $(M,J,\omega)$ be an almost complex manifold of complex dimension 2 with a
%taming symplectic form. Assume that $M$ contains no $J$-holomorphic sphere
%with a self-intersection number equal to $-1$. Suppose that the boundary $\partial M$ of $M$ is a smooth  strictly Levi convex hypersurface. Let $S^2$ be a real $2$-sphere containing  elliptic and 
%hyperbolic points only and embedded into $\partial M$. Then there exists a filling of $S^2$ by boundaries of $J$-holomorphic discs contained in $M$. This filling forms a real Levi-flat hypersurface $\Gamma \subset M$ with boundary $S^2$. 
%\end{theorem} 

\section{Preliminaries}
\subsection{Almost complex and symplectic structures}

All manifolds and almost complex structures are supposed to be of class $C^\infty$
though the main results require a lower regularity.
Let $( M', J')$ and $(M,J)$ be almost complex manifolds and let $f$ be
a map of class $C^1$ from $\tilde M$ to $M$. We say that $f$ is
 $(J',J)$-holomorphic if $df \circ J' = J \circ df$.
 We denote by $\mathbb D$ the unit
disc in $\C$ and by $J_{st}$ the standard  structure on $\C^n$
for every $n$. If $(M',J')=(\mathbb D,J_{st})$,
we call $f$ a  $J$-holomorphic disc in $M$. Every almost complex manifold
$(M,J)$ can be viewed locally as the unit ball $\mathbb B$ in
$\C^n$ equipped with a small almost complex
deformation of $J_{st}$. Indeed for every point $p \in
M$, every real $\alpha \geq 0$ and   $\lambda_0 > 0$ there exist a neighborhood $U$ of $p$ and a
coordinates diffeomorphism $\phi: U \longrightarrow \mathbb B$ such that
$\phi(p) = 0$, $d\phi(p) \circ J(p) \circ d\phi^{-1}(0) = J_{st}$  and the
direct image $ \phi_*(J): = d\phi \circ J \circ d\phi^{-1}$ satisfies $\vert\vert \phi_*(J) - J_{st}
\vert\vert_{C^\alpha(\bar {\mathbb B})} \leq \lambda_0$.

Let $U$ be an open set in $M$ and $\phi : U \longrightarrow \C^n$ be a coordinates diffeomorphism. Then $\phi(U)$ is an open set in $\C^n$ with standard coordinates $z$. In the sequel we often denote $\phi_*(J)$ by $J$. If $f: \mathbb D \longrightarrow U$ is a $J$-holomorphic disc then we still denote by $f$ the composition $\phi \circ f$ viewing $f$ as a $\C^n$-valued map $f :\zeta \mapsto z(\zeta)$. The
conditions $J(z)^2 = -I$ and $J(0) = J_{st}$ imply that for every $z$ the endomorphism $u$ of $\R^{2n}$
defined by  $u(z):= - (J_{st} + J(z))^{-1}(J_{st} - J(z))$ is anti
$\C$-linear that is $ u \circ J_{st} = -J_{st} \circ u$. Thus  $u$ is a
composition of the complex conjugation
and a $\C$-linear operator. Denote by $A_J(z)$ the complex $n \times n$ matrix
such
that $u(z)(v) = A_J(z)  \overline v$ for any $v \in \C^n$.
 The entries of the matrix $A_J(z)$
are smooth functions of $z$ and
\begin{eqnarray}
\label{Norm1}
A_J(0) = 0.
\end{eqnarray} 
The $J$-holomorphicity condition $J(z) \circ df = df \circ J_{st}$ can be
written in the form of the nonlinear Cauchy-Riemann equation~:

\begin{eqnarray}
\label{CR}
\partial_{\bar \zeta} f + A_J(f) \partial_{\bar \zeta} \overline{f}   = 0.
\end{eqnarray}

The matrix function $A_J$ is called the deformation tensor of $J$ in the coordinates $z$.
Consider the isotropic dilations $d_{\lambda}: z \mapsto \lambda^{-1}z$ in $\C^{n}$. Since the structures $J_\lambda:= (d_\lambda)_*(J)$
converge to $J_{st}$ in any $C^\alpha$-norm  as
$\lambda \longrightarrow 0$, we have $\displaystyle A_{J_\lambda} \longrightarrow_{\lambda \longrightarrow 0} 0$ in any
$C^\alpha$ norm. Thus, shrinking $U$ if necessary and using the isotropic
dilations of coordinates we can assume
that for given $\alpha > 0$ we have $\parallel A_J
\parallel_{C^\alpha(\mathbb B)}< <1$ on the unit ball of $\C^n$. In particular, the
system (\ref{CR}) is elliptic.

According to  classsical results \cite{ve}, the Cauchy-Green
 transform

$$Tf(\zeta) = \frac{1}{2\pi i} \int\int_{\mathbb D} \frac{f(\tau)}{\tau - \zeta}d\tau \wedge d\overline\tau$$
is a continuous linear operator from $C^\alpha(\overline{\mathbb D})$ into
 $C^{\alpha+1}(\overline{ \mathbb D})$ for every non-integer $\alpha > 0$ .
Hence  the operator

\begin{eqnarray}
\label{resolution-op}
\Psi_{J}: f \longrightarrow h =  f + TA_J(f) \partial_{\bar \zeta} \overline {f}
\end{eqnarray}
maps the space   $C^{\alpha}(\overline{\mathbb D})$  into itself and we can write   the
equation (\ref{CR}) in the form
$\partial_{\bar \zeta} \Psi_J(f)  = 0$. This means that a disc $f$ is $J$-holomorphic if
and only if the map $h=\Psi_{J}(f)$ is
$J_{st}$-holomorphic.
If the norm of $A_J$ is small enough,
then  by the Implicit Function Theorem the operator $\Psi_J$
realizes a one-to-one
correspondence between sufficiently small $J$-holomorphic discs and
$J_{st}$-holomorphic discs. This implies the existence of a $J$-holomorphic disc
in a given tangent direction through a given point, a smooth dependence of such a
disc  on a deformation of the point, the tangent vector and the almost complex structure, as well as the interior elliptic regularity of
discs. This is the content of the classical Nijenhuis-Woolf theorem, see
\cite{si} for a short complete proof.

 A {\it symplectic form} $\omega$ on a smooth $2n$-dimensional manifold $M$ is a non-degenerate closed differential two-form. The pair $(M,\omega)$ is called a {\it symplectic manifold}. An almost complex structure $J$ on $M$ is called {\it tamed} by $\omega$ if $\omega(v,Jv) > 0$ for every non-zero tangent vector $v$.  Every $\omega$-tamed almost complex structure $J$ defines a Riemannian metric 
$$g_J(u,v) = \frac{1}{2}(\omega(u,Jv) + \omega(v,Ju)).$$

The Gromov Compactness Theorem used in our paper requires some restriction on
$(M,J,\omega)$ called {\it the bounded geometry} condition. This means that
$g_J$ has uniformly bounded sectional curvatures and positive injectivity
radius, uniformly separated from zero, and that $J$ is uniformly continuous with
respect to $g_J$. It is shown in \cite{si} that the bounded geometry condition can be stated in the following form which imposes  less regularity restrictions : 

\bigskip

$(M,J)$ admits a complete Riemannian metric $g$ such that there exist positive constants $r_0$, $C_1$, $C_2$ with the following properties:

\begin{itemize}
\item[(i)] For every $p \in M$ the $\exp_p:B(0,r_0) \to B(p,r_0)$ is a diffeomorphism. Here $B$ denote the corresponding balls.
\item[(ii)] Every loop $\gamma$ in $M$ contained in the ball $B(p,r)$, $r \leq r_0$, bounds a disc in $B(p,r)$ of area less than $C_1{\mathrm  lenght}(\gamma)^2$,
\item[(iii)] On every ball $B(p,r_0)$ there exists a symplectic form $\omega_p$ such that 
$\vert\vert\omega_p\vert\vert \leq 1$ and $\vert X \vert^2 \leq C_2 \omega_p(X,JX)$ (taming property).
\end{itemize}

Properties (i) and (ii) always hold if $M$ is closed and the complete metric $g$ is of class $C^2$ on $M$. Property (iii) holds if $J$ is  additionally uniformly continuous with respect to $g$. In what follows we always assume that the symplectic tamed almost complex  manifold $(M,J,\omega)$ satisfies the bounded geometry condition.

%Everywhere below we assume that $M$ satisfies this condition.

\subsection{Levi form and plurisubharmonic functions} The proof of the statements of this Subsection can be found in \cite{di-su}. 
Let  $r$ be a $C^2$ function on $(M,J)$. We denote by $J^*dr$ the
differential form acting on a vector field $X$ by $J^*dr(X):= dr(JX)$.
For example, if $J = J_{st}$ on $\R^2$, then $J^*dr = r_ydx - r_xdy$.
The value of {\it the  Levi
form of} $r$ at a point $p \in M$ and at a vector $t \in T_p(M)$ is defined  by
$$L^J_r(p;t): = -d(J^* dr)(X,JX)$$ where $X$ is an arbitrary smooth vector field in a
neighborhood of $p$ satisfying $X(p) = t$.  This definition is independent of
the choice of vector fields. For instance, if $J = J_{st}$ in $\R^2$, then
$-d(J^*dr) = \Delta r dx \wedge dy$ ($\Delta$ denotes the Laplacian). In
particular, $L_r^{J_{st}}(0,\frac{\partial}{\partial x}) = \Delta r(0)$.

The following properties of the Levi form are fundamental :

\begin{proposition}
\label{pro1}
Let $r$ be a real function of class $C^2$ in a neighborhood of a point $p \in M$.
\begin{itemize}
\item[(i)] If $F: (M,J) \longrightarrow (M', J')$ is a $(J,J')$-holomorphic map,  and $\varphi$ is a real function of class $C^2$ in a neighborhood of $F(p)$, then for any $t \in T_p(M)$ we have
$L^J_{\varphi \circ F}(p;t) = L^{J'}_\varphi(F(p),dF(p)(t))$.
\item[(ii)] If $f:\mathbb D \longrightarrow M$ is a $J$-holomorphic disc satisfying
  $f(0) = p$, and $df(0)(e_1) = t$ (here $e_1$ denotes the vector
  $\frac{\partial}{\partial \Re \zeta}$ in
  $\R^2$), then $L^J_r(p;t) = \Delta (r \circ f) (0)$.
\end{itemize}
\end{proposition}
 Property (i) expresses the invariance of the Levi form with
respect to biholomorphic maps. Property (ii) is often useful to compute the Levi form if a vector
$t$ is given.

\bigskip

\begin{definition}
\label{psh-def}
A real $C^2$ function $r$ on $(M,J)$ is called
{\it plurisubharmonic} if :
$$
L_r^J(p,t) \geq 0 \mbox{ for any } p \in M \mbox{ and } t \in T_p(M).
$$
A $C^2$ function $r$ is {\it strictly plurisubharmonic} on $M$ if $L_r^J(p,t) > 0$
for any $p \in M$ and $t \in T_p(M) \backslash \{ 0\}$.\smallskip\\
\end{definition}

According to Proposition~\ref{pro1},  $r$ is plurisubharmonic if its composition with any $J$-holomorphic disc is subharmonic.

\medskip
A useful observation due to E.Chirka \cite{ch} is that the Levi form of a function $r$ at
a point $p$ in
an almost complex manifold $(M,J)$ coincides with the Levi form with respect
to the standard structure $J_{st}$ of $\R^{2n}$ if {\it suitable} local
coordinates near $p$ are choosen. See the details in \cite{di-su}. 

\vskip 0,2cm
Let $p$ be a boundary point of a domain $\Omega$ in an almost complex manifold
$(M,J)$. Assume that $\partial \Omega$ is of class $C^2$ in a neighborhood $U$ of $p$.
Then $\Omega \cap U = \{ q \in U: r(q) < 0 \}$ where $r$ is a real function of
class $C^2$ on $U$, $dr(p) \neq 0$.
\begin{definition}
$(i)$ $\Omega$ is called Levi $J$-convex  at $p \in \partial \Omega$
if $L_r^J(p;t) \geq 0$ for any $t \in T_p( \partial \Omega) \cap J(T_p(\partial \Omega))$ and  strictly Levi $J$-convex
 at $p$  if  $L_r^J(p;t) > 0$ for any non-zero $t \in
T_p(\partial \Omega) \cap J(T_p(\partial \Omega))$.  If  $\Omega$ is a relatively compact
domain with $C^2$ boundary in an almost complex manifold $(M,J)$, then $\Omega$
is called Levi $J$-convex  if it is Levi $J$-convex at every boundary
point.
$(ii)$ A real hypersurface $\Gamma = \{ r = 0 \}$ in an almost complex manifold $(M,J)$
is called {\it Levi -flat} if $L_r^J(p,t) = 0$ for every $p \in \Gamma$ and every $t
\in T_p(\Gamma) \cap J(T_p(\Gamma))$.
\end{definition}

This definition does not depend on the choice of defining functions. We simply write a  Levi convex domain dropping $J$ in the notations when an almost complex structure is prescribed.

\bigskip

{\bf Example.} Consider in $(\C^2,J_{st})$ the smoothly bounded domain $\{ (z_1,z_2) \in \C^2: \vert z_1 \vert^2 + \vert z_2 \vert^{2m} < 1 \}$, where $m \geq 1$ is an integer. This domain  is strictly Levi convex for $m = 1$ and Levi convex for $m \geq 1$. Its boundary contains no non-constant $J_{st}$-holomorphic discs for all integer $m \geq 1$.

\bigskip

A strictly Levi convex  domain  admits near every boundary point a strictly plurisubharmonic defining function. One could expect that  Levi convex domains have a similar property i.e. admit local defining plurisubharmonic functions. However, it is well-known that there are smoothly bounded domains in $\C^n$ with Levi convex boundary (for the standard complex structure) which do not admit defining plurisubharmonic functions, see \cite{di-fo}. Fortunately, a weaker property always  holds : a smoothly bounded Levi convex domain in $\C^n$ admits a bounded exhaustion strictly plurisubharmonic function \cite{di-fo}.

The following result obtained in \cite{di-su} is an almost complex analog of the classical result \cite{di-fo}. It is a key for our approach. 

\begin{theorem}
\label{exhaust-theo}
Let $(M,J)$ be an almost complex manifold and let $\Omega
\subset M$ be a relatively compact Levi convex domain  with $C^3$
boundary, such that there exists a $C^2$ strictly plurisubharmonic
function $\psi$  in a neighborhood $U$ of $\partial \Omega$. Let $r$ be any
$C^3$ defining function for $\Omega \cap U$. Then there exist a neighborhood
$U'$ of
$\partial \Omega$ and constants $A > 0$, $0 < \eta_0 < 1$, such that for any $0 <
\eta \leq \eta_0$ the function $\rho = -(-re^{-A\psi})^\eta$ is strictly
plurisubharmonic on $\Omega \cap U'$. If $U$ is a neighborhood of
$\overline \Omega$, then $\rho$ is strictly plurisubharmonic on $\Omega$. Furthermore for a fixed point $p \in \partial \Omega$ and any given $0 < \eta < 1$ there is a neighborhood $U$ of $p$, a strictly
plurisubharmonic function $\psi$ in $U$ and $A > 0$ such that $\rho$ is
strictly plurisubharmonic in $\Omega \cap U$.
\end{theorem}

As a direct consequence we obtain that no $J$-holomorphic disc can touch $\partial \Omega$ from inside.

\begin{definition}
\label{bishop-def}
Let $E$ be a real submanifold in $(M,J)$. A {\it Bishop disc} for $E$ is a $J$-holomorphic disc
$f : \mathbb D \longrightarrow M$ continuous on $\overline{\mathbb D}$ and such
that $f(\partial \mathbb D) \subset E$.
\end{definition}

We will use the following properties of Bishop discs.

\begin{proposition}
\label{translation}
$(i)$ Let $\Omega$ be a smooth Levi convex domain. For any point $p \in \partial \Omega$ there exists a neighborhood $U$ of $p$ with the following property: if $f: \mathbb D \longrightarrow U$ is a Bishop disc for $\partial \Omega$
  then $f(\mathbb D) \subset \overline \Omega \cap U$.

$(ii)$ In the hypothesis of Theorem \ref{exhaust-theo} for any point $p \in \partial \Omega$ there exists a neighborhood $U$ of $p$ with the following property: if $f: \mathbb D \longrightarrow U$ is a Bishop disc for $\partial \Omega$
then either $f(\mathbb D) \subset \partial \Omega$ or the disc $f(\mathbb D)$ is contained in $\Omega$
and is transverse to the holomorphic tangent space of $\partial \Omega$ at every boundary point.
\end{proposition}

\proof Part $(i)$ is proved in \cite{di-su}. Let us prove part $(ii)$. Suppose that $f$ is not contained in $\partial \Omega$. The statement is local so according to Theorem~\ref{exhaust-theo} $\eta$ in the
construction of an exhaustion plurisubharmonic function can be choosen
arbitrarily close to $1$. Fix $\eta = 3/4$. We can assume that in local coordinates $p=0$, the condition~(\ref{Norm1}) is satisfied and $f(1) = 0$. 
Applying the Hopf lemma to the subharmonic function $u(\zeta) = -(-r
e^{-A\psi}\circ f(\zeta))^{3/4}$ on the unit disc we obtain the estimate 
$$\vert r \circ f (\zeta) \vert \geq C (1 - \vert \zeta \vert)^{4/3}.$$ 
We point out that $u$ does not vanish identically since $f$ is not contained
in $\partial \Omega$ so the constant $C$ is strictly positive.
On the other hand the function $r$ is smooth so near the point $1$ we have $r
\circ f(\zeta) =  a (1 - |\zeta|) + O((1 - |\zeta|)^2)$. Hence $a \neq 0$. In particular $f(\partial \D)$ is a real curve in $\partial \Omega$. Let $v$ be a vector tangent to $f(\partial \D)$ at a point $p$. Since $f(\D)$ is a $J$-curve then $Jv$ is tangent to $f(\D)$ at $p$ and does not belong to the real tangent space of $\partial \Omega$ at $p$ by the expression of $r \circ f$. Hence $v$ is transverse to the holomorphic tangent space of $\partial \Omega$ at $p$. Q.E.D.

According to Y.Eliashberg - W.Thurston \cite{el-th} a tangent hyperplane field
$\xi = \{ \alpha = 0 \}$, where $\alpha$ is a 1-form,  on a $(2n+1)$-dimensional manifold $\Gamma$ is
called a {\it positive confoliation} if there exists an almost complex
structure $J$ on the bundle $\xi$ such that
\begin{eqnarray*}
d\alpha(X,JX) \geq 0
\end{eqnarray*}
for any vector $X \in \xi$. The 1-form $\alpha$ is defined up to the
multiplication by a nonvanishing function. Thus, the confoliation condition for $\xi$ is
equivalent to the existence of a compatible Levi convex CR-structure
(in general, non-integrable). In other words, if $\Gamma = \{ r = 0 \}$ is a
smooth Levi convex  hypersurface in
an almost complex manifold $(M,J)$, then the distribution of its holomorphic
tangent spaces $\xi = T\Gamma \cap J(T\Gamma)$ is a confoliation: we can set
$\alpha = J^*dr$. In particular, if $\Gamma$ is a strictly Levi convex
hypersurface, then $\xi = \{ J^*dr = 0 \}$ is a contact structure. Recall that
a tangent hyperplane field $\xi = \{\alpha = 0 \}$ on $\Gamma$ is called a
contact structure if $\alpha \wedge (d\alpha)^n \neq 0$ on $\Gamma$. One of the main
questions
considered by Y.Eliashberg - W.Thurston concerns the possibility to  deform a
given confoliation to
a contact structure or approximate a confoliation by contact structures.
Combining the contact topology techniques with the geometric  foliation theory
 they
obtained several results of this type in the case where $\Gamma$ is of real
dimension 3. The next result of \cite{di-su} works in any dimension.

\begin{theorem}
\label{theo4}
Let $\Omega$ be a relatively compact pseudoconvex domain with $C^\infty$
boundary in an almost complex manifold $(M,J)$. Assume that there exists a
smooth strictly plurisubharmonic function $\psi$ in a neighborhood of $\partial \Omega$.
Then the confoliation of holomorphic tangent spaces $T(\partial \Omega) \cap
J(T(\partial \Omega))$ can be approximated in any $C^{k}$ norm by contact
structures.
\end{theorem}

In a suitable neighborhood of a fixed point $p \in \partial \Omega$
there exists a smooth strictly plurisubharmonic function. So
every  confoliation can be approximated locally by contact structures.

%An important Cancellation Theorem of Eliashberg-Harlamov claims that a sphere with $e$ elliptic and $h$ hyperbolic points is isotopic to a sphere with $e-1$ elliptic and $h-1$ hyperbolic points.

\section{Complex and totally real  points of real surfaces} 
Let $(M,J,\omega)$ be a real four dimensional symplectic manifold with an almost complex structure $J$ tamed by the symplectic form $\omega$. Let $E$ be a real smooth submanifold in $M$. A point $p \in E$ is called {\it totally real} if $T_pE \cap J(p)T_pE = \{0\}$. A submanifold $E$ is called totally real if it is totally real at every point.
Let $S^2$ be a smooth real surface diffeomorphic to the real two-sphere, embedded into $M$.  A point $p \in S^2$ is called
{\it complex} if the tangent space $T_pS^2$ is a $J(p)$-invariant subspace of
$T_pM$. We consider the generic case where $S^2$ admits a finite set $\Sigma$
of complex points so that $S^2 \backslash \Sigma$ is a totally real submanifold
of $M$.
%We will denote the set of totally real points of $SÃ§2$ by $S^2_*$.

Let $p \in S^2$ be a complex point.  We may choose local coordinates centered at $p$ so
that $p = 0$ and the deformation tensor $A_J$ satisfies Condition~(\ref{Norm1}).
The complex $2\times 2$-matrix function $A$
can be written in the form $A_J = [A_J^1,A_J^2]$ where the columns of $A_J$, $A_J^j,\ j=1,2,$ are smooth maps from a neighborhood of the origin in $\C^2$ to $\C^2$. We need the following additional normalization of $A_J$ .

\begin{lemma}
(Adapted coordinates near a complex point). After a suitable  local change of coordinates at $p$ we have : 
\begin{itemize}
\item[(i)] the deformation tensor  $A_J$  satisfies the
normalization condition (\ref{Norm1}) and  
\begin{eqnarray}
\label{Norm2}
A_J^1(z) = L(z_2) + O(\vert z \vert^2)
\end{eqnarray}
where $L:\C \to \C^2$ is an $\R$-linear map,
\item[(ii)] the sphere $S^2$ is given near the origin by the equation
\begin{eqnarray}
\label{ElPoints2}
z_2 = z_1 \overline{z}_1 + \gamma \Re(z_1^2) + o(|z_1|^2) 
\end{eqnarray}
\end{itemize}
\end{lemma}
Such a system of local coordinates is called {\sl adapted}. We use the notation 
$$\rho(z) = z_2 - z_1 \overline{z}_1 - \gamma \Re(z_1^2)  + o(|z_1|^2) $$
We  have used above the standard notation $\phi(x) = o(\psi(x))$  for  the functions $\phi$ and $\psi$ defined in a neighborhood of the origin in $\R$ and  satisfying
$$\lim_{ x \to 0} \frac{\phi(x)}{\psi(x)} = 0.$$
We also use the usual notation $\phi(x) = O(\psi(x))$ if there exists a constant $C > 0$ such that 
$$\vert \phi(x) \vert \leq C \vert \psi(x) \vert$$
in a neighborhood of the origin.

\proof First, choose coordinates so that
$\rho(z) = z_2 + O(z_1^2)$. Keeping this expression we may achieve
Condition (\ref{Norm2}). Then consider a  change of
coordinates $(z_1,z_2) \mapsto (az_1, z_2 + Q(z))$, $a \neq 0$,  with a holomorphic
homogeneous second degree polynomial $Q$ to obtain the expression (\ref{ElPoints2}). It follows from the transformation
rule for the deformation tensor $A_J$ that Condition (\ref{Norm2}) still holds in
the new coordinates. Q.E.D.

%is contained in a strictly pseudoconvex
%hypersurface. Indeed the normalization conditions
%(\ref{Norm1}) and (\ref{Norm2}) imply that a $J$-holomorphic disc $f$ through the
%origin has the following Taylor expansion $f(\zeta) = a \zeta + b\zeta^2 + O(\zeta^3)$
%so that the Levi form of a $C^2$-function at the origin with respect to $J$
%coincides with its Levi form with respect to $J_{st}$. 

The number $\gamma \in [0,+\infty[$, $\gamma \neq 1$, is a local invariant of
$S^2$ : it is independent of a choice of coordinates satisfying Conditions (\ref{Norm1}),
(\ref{Norm2}), (\ref{ElPoints2}).

\begin{definition}
A complex point is called elliptic if $0 \leq \gamma < 1$, parabolic if $\gamma = 1$ and hyperbolic if $\gamma >  1$.
\end{definition}
In the generic case (i.e. after an arbitrary small perturbation ) $S^2$ contains no parabolic point. Thus, in what follows we consider the only case where the set $\Sigma$ of complex points consists of elliptic and hyperbolic points. Denote by $e$ the number of elliptic points and by $h$ the number of hyperbolic points. Since the Euler characteristic of $S^2$ is equal to 2 we have, according to well-known results : $e=h+2$. In particular if $S^2$ contains no hyperbolic point it contains precisely two elliptic points.
%We assume in this Section that the sphere $S^2$ does not contain any hyperbolic point. In particular $S^2$ contains only two elliptic points $p,q \in S^2$ and $S^2 \backslash\{p,q\}$ is the totally real part $S^2_\star$ of $S^2$.

\section{Generation of Bishop discs near an elliptic point}

A local study of Bishop discs near an elliptic point was performed by several
authors \cite{bi,ho,ko-we,su-tu2,ye}. The exposition of this Section is based on the approach developped by A.Sukhov-A.Tumanov in \cite{su-tu2} in any dimension. Here we adapt it to the two-dimensional case.

\subsection{Bishop's equation}
Denote by $\overline \partial_J$ the Gromov operator
$$
\overline \partial_J f:=df + J \circ df \circ J_{st}.
$$

A smooth map $f$
defined on $\D$ is $J$-holomorphic if it satisfies the non-linear Cauchy-Riemann equation $\overline \partial_J f = 0$ on $\D$.

A smooth map $f$
defined on $\D$ and continuous on
$\overline\D$ is a Bishop disc if and only if it satisfies the
following non-linear boundary problem of the
Riemann-Hilbert type for the
quasi-linear operator $\overline\partial_J$:

\begin{displaymath}
(RH): \left\{ \begin{array}{ll}
\overline\partial_J f(\zeta) = 0, \zeta \in \D,\\
\\
\rho(f)(\zeta) = 0, \zeta \in \partial \D.
\end{array} \right.
\end{displaymath}

%\subsection{Bishop's equation near an elliptic point}
First we derive Riemann-Hilbert type boundary problems describing pseudoholomorphic Bishop discs near generic elliptic points
of a real submanifold in an almost complex manifold. 

Since our considerations are local, we suppose that local coordinates are choosen as in Section 1, namely that Conditions (\ref{Norm1}),
(\ref{Norm2}) and (\ref{ElPoints2}) are satisfied near an elliptic point in $S^2$.

Denote by $P(z_1) = z_1 \overline{z}_1 + \gamma \Re(z_1^2) $ the quadratic part of the Taylor expansion of $S^2$ near the origin.
Consider the non-isotropic dilations
$\Lambda_{\delta}: z = (z_1,z_2) \mapsto z' = (\delta^{-1/2}z_1, \delta^{-1}z_2)$. In the new
$z$-variables  (we drop the primes) the image
$S_{\delta}^2 := \Lambda_{\delta}(S^2)$ is
defined by the equation $\rho_\delta (z):=\delta^{-1}\rho((\Lambda_\delta)^{-1} z) = 0$. 
The functions $\rho_\delta$ converge in any $C^k$ norm on every compact
subset of $\C^2$ to the function 
$$\rho^0(Z) =  z_2 -  P(z_1)$$ 
as $\delta \longrightarrow 0$. 
Hence 
the surfaces $S_\delta^2$ converge for the local Hausdorff distance to the model quadric manifold $S_0^2:=\{ \rho^0(Z) = 0\}$, $\C^2$, as $\delta \longrightarrow 0$.

Suppose additionally that $S^2$ is contained in a strictly $J$-Levi convex
hypersurface $\Gamma$ near the origin. Then the Taylor expansion of the
defining function $\phi$ of $\Gamma$ is 
$$\phi(z) = \alpha\Re z_2 + \beta\Im z_2 + O(\vert z \vert^2).$$
Let $f(\zeta) = (z_1(\zeta),z_2(\zeta)) = (a\zeta,0) + O(\zeta^2)$ be a $J$-holomorphic disc tangent to
$\Gamma$ at the origin. Since  $z_2 = O(\zeta^2)$, the $J$-holomorphicity
equations for $f$ in the adapted coordinates imply that there exist $b_1,b_2 \in \C$ such that :
$$
f(\zeta) = (a\zeta+b_1\zeta^2,b_2 \zeta^2) + O(\zeta^3).
$$
%the second order terms in the expansion of $f$ vanish.
Hence the Levi form of $\Gamma$ at the
origin with respect to $J$ coincides with the Levi form of $\Gamma$ with respect to
$J_{st}$. In particular :
$$
\frac{\partial^2\phi}{\partial z_1 \partial \overline z_1}(0) > 0.
$$
The dilated hypersurface $\Gamma_\delta := \{
\rho_\delta(z) = 0 \}$ converges to a real quadric hypersurface $\Gamma_0$ which is $J_{st}$-strictly
pseudoconvex near the origin.

Consider the pushed-forward structures $J_\delta := (\Lambda_\delta)_*(J) = d\Lambda_\delta \circ J \circ (d\Lambda_\delta)^{-1}$.

\begin{lemma}
(Structure deformation near a complex point) \begin{itemize}
\item[(i)] For every positive integer $k$ and every compact subset $K \subset \C^n$ we have $\parallel J_\delta - J_{st} \parallel_{C^k(K)} \longrightarrow 0$ 
as $\delta \longrightarrow 0$. Thus $S^2$ and $J$ are small deformations of
$S^2_0$ and $J_{st}$, respectively, near the origin.
\item[(ii)] Suppose that $S^2$ is contained in a strictly $J$-Levi convex hypersurface $\Gamma$ near the origin. Then $\Gamma$ is a small deformation of a strictly $J_{st}$-pseudoconvex hypersurface $\Gamma_0$ containing
  $S^2_0$. 
\end{itemize}
\end{lemma} 
\proof (i) Consider the Taylor expansion of $J(z)$ near the origin: $J(z)
= J_{st} + L(z) + R(z)$ where $L(z)$ is  the linear part of the expansion
 and $R(z)= O(\vert z \vert^2)$. Clearly, $\Lambda_{\delta} \circ
R(\Lambda_{\delta}^{-1}(z)) \circ \Lambda_{\delta}^{-1}$ converges to $0$
as $\delta \longrightarrow 0$. Fix $j,k \in \{1,2\}$ and denote by $L_{kj}^{\delta}(z)$
(respectively, by $L_{kj}(z)$) an entry
of the real  matrix $\Lambda_{\delta} \circ
L(\Lambda_{\delta}^{-1}(z)) \circ \Lambda_{\delta}^{-1}$
(respectively, of $L(z)$). We have $J_{12}^{\delta}(z) = \delta^{1/2}L_{12}(\delta^{1/2}
z_1,\delta z_2)$ and  $J_{jj}^{\delta}(z) = L_{jj}(\delta^{1/2}
z_1,\delta z_2)$, $j=1,2$. Hence, these terms tend to $0$ as   $\delta
\longrightarrow 0$. Furthermore  $J_{21}^{\delta}(z) = \delta^{-1/2}L_{21}(\delta^{1/2}
z_1,\delta z_2) \longrightarrow L_{21}(z_1,0)$ as $\delta \longrightarrow
0$, uniformly on $K$. However, it follows from (\ref{Norm2}) that  
$L_{kj}(z_1,0) = 0$ for $k,j =1,2$. This implies that 
$L_{kj}^{\delta}$ converges to $0$ as $\delta \longrightarrow 0$, for every $k,j=1,2$. Now Part (ii) follows from Part (i). Q.E.D.

Let $f$ be a $J_\delta$-holomorphic disc
in a neighborhood of the origin in $\C^2$. The
$J_\delta$-holomorphicity condition for $f$ has the form (\ref{CR}) with the deformation tensor $A_{J_\delta}$ associated to $J_\delta$.

Considering the operator $\Psi_{J_\delta}$ defined by (\ref{resolution-op})
we can replace the non-linear Riemann-Hilbert problem
(RH) by  the Bishop equation

\begin{eqnarray}
\label{Bishop2}
\rho_\delta(\Psi_{J_\delta}^{-1}(h))(\zeta) = 0, \zeta \in \partial \D
\end{eqnarray}
for an unknown function $h$, holomorphic in $\D$ and continuous on $\overline{\D}$.

If $h$ is a solution of the boundary  problem (\ref{Bishop2}) then $f := \Psi_{J_\delta}^{-1}(h)$ is a Bishop disc with  boundary
attached to $S^2_\delta$. Since the manifold $S_\delta^2$ is
biholomorhic via  non-isotropic dilations to the initial manifold $S^2$,
the solutions of Equation (\ref{Bishop2}) allow to describe
the Bishop discs attached to the sphere  $S^2$ near an elliptic point.

\subsection{Geometry of the Bishop discs near an elliptic point}

Our goal is to  prove the following statement :

\begin{theorem}
\label{ellipticfilling-thm}
\begin{itemize} 
\item[(i)] Let $p$ be an elliptic point in $S^2$. Given a positive integer $k$
  and $0 < \alpha < 1$ there exists a family $(f^t)_t$
of  $J$-holomorphic Bishop discs for $S^2$, $C^{k,\alpha}$-smoothly depending on one real parameter $t \in [0,t_0]$. These discs foliate a real hypersurface
$E$ such that $(E, S^2)$ is a $C^{k,\alpha}$ smooth manifold with boundary, outside $p$.
\item[(ii)] Suppose additionally that $S^2$ is contained in the boundary
$\partial\Omega$ of a Levi $J$-convex domain. Then the generated family of discs
  is contained in $\overline\Omega$.
\end{itemize}
\end{theorem} 
Part (ii) of Theorem~\ref{ellipticfilling-thm} follows from part $(i)$ of Proposition~\ref{translation} so we prove the part (i).

We begin with the description of the Bishop discs attached to the model quadric manifold $S_0^2$ in $\C^2$ with the standard structure $J_{st}$.
They are the solutions of the boundary problem (\ref{Bishop2}) for $\delta = 0$.

For $r > 0$  consider the ellipse $D_r: = \{ \zeta \in \C: P(\zeta) < r \}$  and denote by $z_{1,r}$ the biholomorphism 
$z_{1,r} = z_1(r,\bullet):\D \longrightarrow D_r$ satisfying $z_{1,r}(0) = 0$, $(\partial z_{1,r} / \partial \zeta)(0)
> 0$. Then $P \circ z_{1,r} \vert_{ \partial \D}  \equiv r$ and we set $z_2(r,\zeta)
\equiv r$. The maps $\zeta \mapsto (z_1(r,\zeta), z_2(r,\zeta))$ provide a one parameter family of $J_{st}$-holomorphic Bishop discs.
Their boundaries are disjoint  and fill a pointed neighborhood of the origin in $S_0^2 \backslash \{ 0 \}$. For $r = 0$
these discs degenerate to the constant map $\zeta \mapsto (0,0)$. Their images form a foliation of a real  hypersurface  $E$  
such that $(E,S_0^2)$ is a smooth manifold with boundary outside the origin. 

We claim that in the general case of an almost complex structure $J$ the Bishop discs have similar properties. Indeed, for a sufficiently small real positive number $\delta$ let $p^\delta:=(0, \delta)$ be a point on the real ``normal'' to $S^2$. The image of $p^\delta$ by the non-isotropic dilation $\Lambda_\delta$ coincides with 
$p^0:= (0,1)$. There exists a unique Bishop disc $f^0=(z_1^0,z_2^0)$ of the described above family centered at $p^0$; it corresponds to the parameter $r = 1$. 
The parametrizing map $F:(r,\zeta) \longrightarrow (z_1(r,\zeta),z_2(r,\zeta))$ has maximal rank when the parameter $r$ 
is in a neighborhood $U$  of the point $1$ and $\zeta \in \D$. Linearize the equation (\ref{Bishop2}) for $\delta = 0$, at the disc
$f^0$. We claim that the linearized operator is surjective between
the corresponding Banach spaces. Indeed, the linearized Bishop equation
$$
\Re
(\partial_z \rho_0(f^0)\dot z) = \psi
$$
has the form

$$
\left\{
\begin{array}{lll}
-\Re (z_1^0 -  \gamma z_1^0)(e^{i\theta})\dot z_1(e^{i\theta}) + (1/2)\Re \dot z_2(e^{i\theta}) & = &
\psi_1(e^{i\theta}),\\
& & \\
\Re (1/2i)\dot z_2(e^{i\theta}) & = & \psi_2(e^{i\theta}).
\end{array}
\right.
$$

The second equation admits a one-parameter solution given by the Schwarz
integral
\begin{eqnarray*}
(1/2i)\dot z_2(\zeta) = \frac{1}{2\pi i}\int_{\partial\D} \frac{\tau +
  \zeta}{\tau - \zeta}\frac{\psi_2(\tau)}{\tau}d\tau + ic_1
\end{eqnarray*} 
with $c_1 \in \R$. Since every ellipse $D_r$ is homotopic to the
unit disc, the winding number of the function $\zeta \mapsto  (z_1^0 - 
\gamma z_1^0)$ is equal to the winding number of $z_1^0(\zeta) \equiv \zeta$ i.e. is
equal to $-1$. Solving then the boundary problem for $\dot z_1$ we obtain a general
solution depending on four real parameters. If $\psi_j$ are of class
$C^{k,\alpha}(\partial\D)$ for some positive integer $k$ and $0 < \alpha < 1$,
then $\dot z \in C^{k,\alpha}(\overline\D)$ by the classical regularity
properties of the Cauchy type integral \cite{ve}. Applying the Implicit Function Theorem to the operator equation (\ref{Bishop2}) we obtain for every sufficiently small positive real number $\delta$ 
a $J_{st}$-holomorphic solution $h(r,\bullet)$, smoothly depending on four real parameters. Three of these four parameters may be removed by fixing a parametrization of the discs. Thus, we obtain a
one-parameter family of disjoint discs. Then $f(r) = \Psi_{J_\delta}^{-1}(h)$ is a family of $J$-holomorphic Bishop discs with boundaries attached to $S_\delta^2$. This family is a small deformation of the above Bishop discs attached to $S_0^2$. Since $\delta$ is small enough, the parametrizing map 
$(r,\zeta) \mapsto f(r,\zeta)$ has maximal rank. So these discs swept a real smooth hypersurface $\Sigma_\delta$ with smooth boundary.
In order to conclude the proof of theorem  it remains  to show  that  
this hypersurface is foliated by the Bishop discs $f(r,\bullet)$. This is a consequence of the general  uniqueness result which we establish in the next subsection. It will also have further applications in our approach.

\subsection{Uniqueness of Bishop discs} 
Uniqueness results for Bishop discs were obtained by several authors in different forms (see for instance \cite{be-ga, be-kl, ye}). Here we follow the exposition of \cite{ye}.
 
\begin{proposition}
\label{ellipticuniqueness-prop}
Let $\Omega$ be a Levi convex domain, relatively compact in $(M,J)$.
Let $S^2$ be an embedded sphere in $\partial \Omega$ and $p$ an elliptic point in $S^2$. Suppose that $\partial \Omega$ contains no germ of a nonconstant $J$-holomorphic disc near $p$. Then there exists a neighborhood $U$ of $p$ in $M$ with the following property :
if $f:\D \to  \Omega \cap U$ is a $J$-holomorphic disc with boundary attached to $S^2 \backslash \{ p \}$, then the image of $f$ coincides with one of the discs of the family constructed in Theorem \ref{ellipticfilling-thm}.
\end{proposition}
\proof We proceed in several steps.

{\it Step 1. Transversality to the boundary.} Since $\partial \Omega$ is Levi convex, it follows from the Hopf lemma that every $J$-holomorphic disc in $\Omega$ with boundary glued to $\partial \Omega$ intersects $\partial \Omega$ transversally at all points of the boundary. 

\vskip 0,2cm
{\it Step 2. Stability of intersections}. This is the content of the following 

\begin{lemma} 
\label{ElPointsStability}
Let $f$ and $g$ be two distinct $J$-holomorphic Bishop discs glued to $S^2 \backslash\{p\}$ such that they intersect at a boundary point $q \in S^2$. If $\tilde f$ is a Bishop disc close enough to $f$ in the $C^2$-norm, then $\tilde f(\overline\D)$ has a non-empty intersection with $g(\overline\D)$. The same holds if they intersect at an interior point where at least one of them is immersed.
\end{lemma}
\proof If $f$ and $g$ intersect transversally at $q$ then their boundaries at $q$ are transverse and the assertion is obvious. Assume that they are tangent at $q$. We can suppose  after a reparametrization of the maps $f$ and $g$ that they are defined on the half-disc $\D^+ = \{ \zeta: \Im \zeta > 0 \}$ and are smooth up to the boundary.  Then we can choose local coordinates near $q$ such that $q=0$, $J(0) = J_{st}$, $S^2 = \{z \in \C^2: \Im z = 0\}$ and $g(\zeta) = (\zeta,0)$. Then $f$ and $\tilde f$ are the graphs over $g(\D^+) = \D^+ \times \{ 0 \}$ of the functions $h$ and $\tilde h$ respectively : $f(\zeta)=(\zeta,h(\zeta))$ and $\tilde f(\zeta)=(\zeta,\tilde h(\zeta))$ for $\zeta \in \D^+$. The tangency condition implies 
\begin{eqnarray}
\label{ElPoints12}
h(\zeta) = a \zeta^k + o(\zeta^k)
\end{eqnarray}
for some  integer $k \geq 2$, with $a \in \R$. Since the boundary of $f$ is glued to $S^2$, we have $\Im h\vert_{[-1,1]} \equiv \Im \tilde h\vert_{[-1,1]} \equiv 0$. Then we may extend $h$ and $\tilde h$ continuously to $\D$ setting $h(\zeta) = \overline h(\overline \zeta)$ for $\zeta \in \D \backslash \D^+$. The extension of $h$ still satisfies (\ref{ElPoints12}) and by the mapping degree theory $\tilde h$ has exactly $k$ zeros, counted with their multiplicities, near the origin. But since it is obtained by reflection over the real axis, it admits at least one zero in $\overline {\D^+}$.
The proof for an interior point is similar with obvious simplifications. Q.E.D.

\vskip 0,2cm
{\it Step 3.} We prove Proposition~\ref{ellipticuniqueness-prop}. Consider the one parameter family $(f^t)_t$ of Bihop discs constructed in Theorem~\ref{ellipticfilling-thm}. For $t$ sufficiently close to $0$ the Bishop disc $f^t$ does not intersect $f$. Let $\tau > 0$ be the infimum of $t \in [0,t_0]$
such that $f$ intersects $f^t$. Then by continuity $f$ intersects $f^\tau$ at an interior point or at a boundary point (notice that the discs $f^t$ are embedded). If $f$ and $f^\tau$ are distinct then Lemma \ref{ElPointsStability} gives a contradiction to the minimality of $\tau$. Q.E.D.

\section{Deformation of discs with totally real boundaries} The results of
this Section can be deduced from general Theorems due to H.Hofer \cite{ho},
H.Hofer-V.Lizan-J.C.Sikorav \cite{ho-li-si} and R.Ye \cite{ye}. Since our situation is rather special we give a direct approach.

Consider a $J$-holomorphic embedding $f^0: \D \to (M,J)$ of class $C^{k,\alpha}(\overline \D)$ for some $k \geq 1$ and $0 < \alpha < 1$. We suppose that $f^0(\partial \D)$ is contained in a totally real submanifold $E$ of dimension $2$ in $M$.  Furthermore we assume that $E$ is contained in a smooth real hypersurface $\Gamma$. We always assume that $f^0$ is transverse to $\Gamma$. The aim of this Section is to generate Bishop's discs attached to $E$ by deforming $f^0$. This will be a direct application of the Implicit Function Theorem. We proceed in several steps.

%Below we will introduce a homotopic invariant of $f^0$ called {\it the winding number} of $f$. The aim of this Subsection is to prove the following
 
%\begin{proposition}
%\label{discdeformation-prop}
%Suppose that the winding number of $f^0$ is equal to $0$. Then there exists a one parameter family $(f^t)_t$ of Bishop discs for $E$, close to $f^0$ in the $C^{k,\alpha}(\overline{\D})$-norm, such that the boundaries of $(f^t)_t$ fill a neighborhood of $f^0(\partial \D)$ in $E$.
%\end{proposition}

%Proposition~\ref{discdeformation-prop} is a direct application of the Implicit Function Theorem. We proceed in several steps.
 
{\it Step 1. Choice of coordinates.} Let $h_\zeta:\D \to M$ be a family of $J$-holomorphic discs, smoothly depending on a parameter $\zeta \in \partial \D$, such that $h_\zeta(0) = f(\zeta)$ for every $\zeta \in \partial \D$ and $h_\zeta$ is tangent to $\Gamma$ at the point $f(\zeta)$. Since $f^0$ is an embedding there exists a neighborhood $U$ of $f^0(\overline\D)$ and a coordinate diffeomorphism $H$ from $U$ to $H(U) \subset \C^2$  such that

$$
(H \circ f^0)(\zeta) = (\zeta,0), \zeta \in \D
$$

and 

$$
J \vert_{H \circ f^0(\D)} = J_{st}
$$

i.e. the deformation tensor $A_J$ of $J$ satisfies

\begin{eqnarray}
\label{normalisation}
A_J(\zeta,0) = \partial_{z_1}A_J(\zeta,0) = 0
\end{eqnarray}
(see \cite{su-tu1}). Furthermore for every $\zeta \in \partial \D$ the map $\phi \circ h_\zeta : \tau \in \D \mapsto (H \circ h_\zeta)(\tau)$ has the form $(H \circ h_\zeta)(\tau) = (\zeta,\tau)$.   We call these coordinates {\it normal coordinates} along
the disc $f^0$.
\vskip 1cm

For simplicity of notations we still denote by $f^0$ the composition $H \circ f^0$ and by $h_\zeta$ the composition $H \circ h_\zeta$.

\vskip 0,2cm
{\it Step 2. Winding number.} The manifold $E$, in a neighborhood of the circle $\partial \D \times \{ 0 \}$, forms a bundle over this circle with as fibers smooth curves $\gamma_\zeta$ which are tangent to the the complex plane $\{ \zeta \} \times \C$  at $(\zeta,0)$, $\zeta = e^{i\theta} \in \partial \D$. Denote by $X_1(\zeta)$ the vector field $(i\zeta,0)$ tangent to the circle $\partial \D \times \{ 0 \}$. Let also $X^2(\zeta)$, $\zeta \in \partial \D$, be a vector field tangent to the fiber $\gamma_\zeta$ at the point $(\zeta,0)$. We assume that $E$ is orientable along $f^0(\partial \D)$; this assumption will always hold in our applications.
Then $X^2(\zeta) = (0,ie^{i\phi(\theta)})$ where $\phi$ is a
smooth function on $[0,2\pi]$ such that $\phi(0) = \phi(2\pi)$. Let $\mu:=\mu(f^0)$ be the winding number of the function
$\theta \mapsto e^{i\phi(\theta)}$, or equivalently the number of times which the tangent vector to $\gamma_\zeta$
turns around the origin in the plane $\{ \zeta \} \times \C$ when $\zeta$ runs
over $\partial \D$. Denote by $L_\zeta$ the real tangent space of $E$ at
$f(\zeta)$, generated by $X^1(\zeta)$ and $X^2(\zeta)$.
%Thus we can view $L_\zeta$ as a vector bundle over the unit circle $\partial \D$.

It is worth pointing out that $\mu$ is defined independently of a choice of normal coordinates. Obviously $\mu$  is invariant with respect to the homotopy of Bishop discs. Namely, if $f^t$ is a family of Bishop discs for $E$ continuously depending on a real parameter $t$, then  normal coordinates corresponding to $f^t$ can be choosen smoothly depending in $t$. Therefore $\mu(f^t)$ continuously depends in $t$ and so is a constant function in  $t$.

\vskip 0,2cm
Let us return to the family of Bishop discs constructed near an elliptic point.
%in the previous section.

\begin{lemma}
\label{maslov-lemma}
The winding number $\mu(f)$ of a Bishop disc $f$ near an elliptic point is equal to $0$. 
\end{lemma}

\proof By the non-isotropic rescaling every Bishop disc is homotopic to the
disc $(z^1_r(\zeta),r)$ glued to the model quadric $S_0^2$. Its winding number 
is equal to $0$ which implies the statement. Q.E.D.

Hence, for all $t$ we have $\mu(f^t) = 0$.

%Indeed, let $H_\zeta(\Gamma)$ denotes the holomorphic tangent space of $\Gamma$ at $f(\zeta)$ i.e. $H_\zeta(\Gamma) =
%T_{f(\zeta)}(\Gamma) \cap J(f(\zeta))T_{f(\zeta)}$. Since $H_\zeta(\Gamma)$ is a complex space, it has the natural orientation.
%Consider the line $l(\zeta) = T_{f(\zeta)}E \cap H_\zeta(\Gamma)$. Then $k$
%counts the number of rotations (in accordance with the   the orientation) of the line   $l(\zeta)$  around the origin in the oriented space $H_\zeta(\Gamma)$. Then number $k$ is called {\it the Maslov index} of the $E$ along the loop $f^0(\partial \D)$. Obviously, it is invariant under the homotopies of the map $f^0$.

\vskip 0,2cm
{\it Step 3. Linearization.}  Now we are able to prove the main result of this section.

\begin{proposition}
\label{discdeformation-prop}
Suppose that the winding number of $f^0$ is equal to $0$. Then there exists a one parameter family $(f^t)_t$ of Bishop discs for $E$, close to $f^0$ in the $C^{k,\alpha}(\overline{\D})$-norm, such that the boundaries of $(f^t)_t$ fill a neighborhood of $f^0(\partial \D)$ in $E$.
\end{proposition}

 If $f: \zeta \mapsto (z_1(\zeta),z_2(\zeta))$ is a $J$-holomorphic disc close to $f^0$, the Cauchy-Riemann equations satisfied by $f$ have the form

$$
\left\{
  \begin{array}{lll}
    \partial_{\bar \zeta} z_1 - a_{11}(z) \partial_{\bar \zeta} {\overline z}_1 - a_{12}(z) \partial_{\bar \zeta} {\overline z}_2 & = & 0, \quad \zeta \in \D\\
    & & \\
    \partial_{\bar \zeta} z_2 - a_{21}(z) \partial_{\bar \zeta} {\overline z}_1 - a_{22}(z) \partial_{\bar \zeta} {\overline z}_2 & = & 0, \quad \zeta \in \D.
\end{array}
\right.
$$

Therefore in order to prove that the linearized Cauchy-Riemann operator is surjective we must solve the following system of linear PDEs  :

$$
\left\{
  \begin{array}{lll}
    \partial_{\bar \zeta} \dot z_1 - b_{11}(\zeta) {\dot z}_2 - b_{12}(\zeta) \overline{\dot z}_2 & = & h_1, \quad \zeta \in \D\\
    & & \\
    \partial_{\bar \zeta} \dot z_2 - b_{21}(\zeta) {\dot z}_2 - b_{22}(\zeta) \overline {\dot z}_2 & = & h_2, \quad \zeta \in \D.
  \end{array}
\right.
$$
Here $h=(h_1,h_2)$ is a given map of class  $C^{k-1,\alpha}(\D)$, the coefficients $b_{ij}$ are of class $C^{k,\alpha}(\D)$ and $\dot z:\D \to \C^2$ denotes an unknown map of class $C^{k,\alpha}(\D)$.
The above special form of the linearized equations along $f^0$ is a
consequence of  the normalization condition (\ref{normalisation}).

Moreover if $E$ is defined by the equation $\rho := (\rho_1,\rho_2) = 0$,
the linearization of the boundary condition $\rho \circ z \vert_{\partial \D} = 0$ has the form 
$$
\dot z(\zeta) \in L_\zeta, \quad \zeta \in \partial \D.
$$
This is equivalent to the following conditions for $\theta \in [0,2\pi]$ :
$$
\Re (e^{-i\theta} \dot z_1(e^{i\theta})) = 0,\  \Re (e^{-i\phi(\theta)} \dot z_2(e^{i\theta})) = 0.
$$

In particular the non-homogeneous boundary conditions are for $\theta \in [0,2\pi]$ :

\begin{eqnarray}
& &\Re (e^{-i\theta} \dot z_1(e^{i\theta})) = g_1(e^{i\theta}),\\
& &\Re (e^{-i\phi(\theta)} \dot z_2(e^{i\phi(\theta)})) = g_2(e^{i\theta}).
\end{eqnarray}

The linearized boundary value problem splits into two parts. Consider first the problem
for $\dot z_2$.
Since $\mu = 0$, It follows from classical results (see \cite{ve}) that the boundary value
problem
$$
\left\{
  \begin{array}{lll}
     \partial_{\bar \zeta} \dot z_2 - b_{21}(\zeta) {\dot z}_2 - b_{22}(\zeta) \overline {\dot z}_2 & = & h_2, \quad \zeta \in \D\\
          & & \\
     \Re (e^{-i\phi(\theta)} \dot z_2(e^{i\theta})) & = & g_2(e^{i\theta})
  \end{array}
\right.     
$$
is  solvable in the class $C^{k,\alpha}(\overline \D)$ for any $h_2 \in C^{k-1,\alpha}(\D)$ and $g_2 \in
C^{k,\alpha}(\partial \D)$. Furthermore the corresponding 
homogeneous problem admits precisely  one  linearly independent solution. Considering then $\dot z_2$ as a
known function, the boundary problem for $\dot z_1$
admits a general solution of class $C^{k,\alpha}$ depending on 3 real
parameters. Thus the linearized non-homogeneous problem always admits a
solution. By the Implicit
Function Theorem there exists a family of discs, depending on
$4$ real parameters, glued to $E$ in a neighborhood of $f^0$. Three parameters correspond to a parametrization of $\D$ and must be removed if we seek for discs with different images. In particular if the winding number $\mu$ is equal to zero we obtain a one-dimensional family of discs forming a one dimensional submanifold ${\mathcal M}$ in the Banach space $C^{k,\alpha}(\D)$. Consider the evaluation map $ev: {\mathcal M} \to E$ 
defined by $ev(z) = z(1)$. We claim that it has maximal rank equal to 1 at $z = (\zeta,0)$. Indeed, the tangent map $\dot ev$ at the disc $f^0=(\zeta,0)$ is defined on $T_{f^0}(\mathcal M)$ by $\dot ev : \dot z \mapsto \dot z(1)$. In order to give a suitable parametrization of the tangent space $T_{f^0}(\mathcal M)$ consider the following auxiliary problem :%to ${\mathcal M}$ at $f^0$ is formed by the maps $\dot z$ satisfying the linearized Cauchy-Riemann equations with homogeneous boundary conditions :

$$
\left\{
\begin{array}{lll}
\partial_{\bar \zeta} w - b_{22}(\zeta) w - b_{22}(\zeta) \overline w = h, \quad \zeta \in \D & & \\
& & \\
\Re   w(e^{i\theta}) = 0 & &
\end{array}
\right.
$$
where $h$ is a given function of class $C^{k-1,\alpha}(\overline \D)$ and $w$ is an unknown function of class $C^{k,\alpha}(\overline{\D})$.
\vskip 0,1cm
By the classical results of Bojarski \cite{bo} that boundary problem admits a unique solution satisfying the condition $\Im w(1) = 0$. Then given $c \in \R$ set $h =   b_{22}(\zeta) ic - b_{22}(\zeta) ic$. If $w$ is a solution of the boundary problem with $\Im w(0) = 0$, then $\dot z_2 = w + ic$ satisfies the linearized homogeneous problem and $\dot z_2(0) = ic$. This implies that the evaluation map has rank 1. Therefore the boundaries of the discs in ${\mathcal M}$ fill a neighborhood of the circle $f^0(\partial \D)$ on $E$.
%This proves Proposition~\ref{discdeformation-prop}. Q.E.D.

%\subsection{Maslov index and the winding number} In conclusion we recall the
%notion of the Maslov index and its relation to the winding number. Denote by
%${\Theta}_2 = GL(2,\C)/Gl(2,\R)$ the manifold of totally real 2-dimensional
%subspaces of $\C^2$. Denote also by $L_0$ the loop 
%$$L_0(e^{2\pi i\theta}) = \{ e^{\pi i \theta} \R \oplus \R \}, \theta \in [0,1]$$
%Then the map 
%$$\mu_0: \pi_1({\Theta}_2) \to {\mathbb Z}$$
%$$\mu_0(L_0) = 1$$
%is an isomorphism. One can define it explicitely. Namely, consider the manifold ${\Lambda}_2$ of Lagrangian $2$-dimensional subspaces of $\C^2$. Consider the map  
%$\pi:{\Theta}_2 \to S^1$, $S^1$ being the unit circle, defined by $\pi( a \cdot GL(2,\R) = \det (a^2)/\det (a^* a)$ where $a^*$ denotes the conjugate transpose. If $L = (L_z)$, $z \in S^1$ is a loop in ${\Theta}_2$ then $\mu_0(L) = deg (\pi \circ L)$. The number $\mu(L)$ is called the {\it Maslov index} of a loop $L$.

%Let now $f^0$ be an embedded  $J$-holomorphic Bishop disc with the boundary
%contained in the totally real part of $S^2$. Let  $L$ be the loop defined by the tangent spaces of $S^2$ along the boundary of $f^0$. In this case the Maslov index of $L$ is called {\it the Maslov index of the Bishop disc $f^0$}.
%One can show (see \cite{ho}) that 
%$$\mu_0(L) = \mu(f^0) + 2$$
%i.e. in our case the Maslov index of Bishop discs generated by elliptic points is equal to $2$.

\section{Filling of a sphere with two elliptic points}
This Section is devoted to the proof of Theorem~\ref{elliptic-thm}.
Since the proof relies on the Gromov Compactness Theorem we recall related notions following \cite{mc-sa, si}. Let $S$ be a compact Riemann surface with (possibly empty) smooth boundary $\partial S$.
We use the canonical identification of the
complex plane $\C$ with ${\bf CP} \backslash \{ \infty \}$. Let
$(M,J,\omega)$ be a symplectic manifold with a tamed almost complex
structure. We assume that $M$ has bounded geometry. Let $E$ be a smooth compact totally
real submanifold of maximal dimension in $M$.

Consider a sequence $f_n: S \to M$ of $J$-holomorphic maps.

Let $g: {\bf CP} \to M$ be a non-constant $J$-holomorphic map. We
say that $g$ occurs as a sphere bubble for the sequence $(f_n)$ if there exists a sequence of holomorphic charts $\phi_n :R_n\D \to S$ with $R_n \to \infty$
converging uniformly on compacts subsets of $\C$ to a constant map $\phi_\infty = p \in S$ and such that  
$$f_n \circ \phi_n \to g$$
uniformly on compact subsets of $\C$.

Let $g:\D \to M$ be a $J$-holomorphic map, continuous on
$\overline\D$, with $\partial \D \subset E$. We say that $g$ occurs as a disc bubble for the sequence $(f_n)_n$ if
there exists a sequence of holomorphic charts $\phi_n:\D \backslash (-1 + \delta_n\D) \to
  M$, smooth on $\overline\D \backslash ( -1 + \delta_n\D )$ with 
$\phi_n (\partial\D   \backslash ( -1 + \delta_n\D )) \subset E$ and
$\delta_n \to 0$, such that $(\phi_n)$ converge uniformly on compact subsets
of $\overline \D \backslash \{ - 1 \}$ to a constant map point $\phi_\infty = p \in S \cup \partial S$ and 
$$f_n \circ \phi_n \to g$$
uniformly on compact subsets of $\overline \D \backslash \{ - 1 \}$.

We have the following simple version of Gromov's Compactness
Theorem :

\begin{proposition}
\label{Gromov1}
Let $f_n:\D \to M$ be a sequence of $J$-holomorphic diss continuous on $\bar\D$
satisfying the boundary conditions  $f_n(\partial \D) \subset E$, intersecting a fixed compact
subset $ K \in M$ and such that 
$$area (f_n):=\int_{\D} f_n^* \omega \leq c$$
where $c > 0$ is a constant. Then there exists a finite set $\Sigma$ in $\bar\D$, eventually empty, such that after extraction :
\begin{itemize}
\item[(i)] The sequence $(f_n)_n$ converges uniformly on compact subsets of $\bar\D \backslash \Sigma$ to a
  $J$-holomorphic map $f_\infty:\D \to M$,
\item[(ii)] A bubbling of some $J$-holomorphic sphere occurs at every point in $\Sigma \cap \D$,
\item[(iii)] A bubbling of some $J$-holomorphic disc occurs at every point in $\Sigma \cap \partial \D$.
\end{itemize}
\end{proposition}

\noindent{\bf Proof of Theorem \ref{elliptic-thm}}.
We proceed in several steps.

\vskip 0,1cm
{\it Step 1. Characteristic foliation.} Let $p$ and $q$ be the elliptic points in $S^2$ and let $f^t$, $g^t$ be the families of Bishop
discs near these points, given by Theorem~\ref{ellipticfilling-thm}. Denote by $\chi$ the bundle over $S^2$ formed by
$T(\partial M) \cap J T(\partial M)\cap T(S^2)$.  The integral
curves of this bundle form a foliation (the characteristic foliation) of $S^2
\backslash \{ p,q \}$ with two singular points $p$ and $q$. Every Bishop disc
is transverse to $\chi$ by part $(ii)$ of Proposition~\ref{translation}. In particular, the characteristic foliation cannot have any closed trajectories  and the boundary of a Bishop disc intersects every
leaf of the characteristic foliation at a single point.  Fix such an integral
curve parametrized by a parameter $t \in [0,1]$. We may assume that the
family of Bishop discs $(f^t)$ near $p$ is parametrized by $t$ close to
$0$.

\vskip 0,1cm
{\it Step 2. Area estimate.} Let $T$ be the supremum of parameters $t$ for which the family $(f^t)_t$ is defined. The areas of the discs $f^t,\ t \in [0,T[$, are uniformly bounded with respect to $t$. Indeed for a fixed $t$ the
boundary of $f_t$ divides $S^2$ into two discs, $S^2_+(t)$ and $S^2_-(t)$. Since
$\omega$ is closed, the Stokes Theorem implies 
\begin{eqnarray*}
\int_{f_t(\D)} \omega = \int_{S^2_+(t)} \vert\omega \vert \leq   \int_{S^2} \vert
\omega \vert.
\end{eqnarray*}

Hence it follows from Proposition~\ref{Gromov1} that there is a sequence $t_k \longrightarrow T$ as $k \longrightarrow \infty$ and a finite set $\Sigma \subset \overline \D$ such that the sequence $(f_k) = (f^{t_k})_k$ converges uniformly to a $J$-holomorphic disc on every compact subset of $\overline \D \backslash \Sigma$. Moreover if $\Sigma \neq \emptyset$ then every point in $\Sigma$ corresponds to a disc bubble or to a sphere bubble.
Since $M$ does not contain any holomorphic sphere by assumption then $\Sigma \cap \D = \emptyset$. 

\vskip 0,1cm
{\it Step 3. Non appearance of bubbles.}
We have the following

\begin{lemma}
\label{nobubble-lem}
There is a sequence of Mobius transformations $\phi_k:\D \to \D$ such that (after extraction) the discs $f_k \circ \phi_k$ converge uniformly on $\bar\D$ to a $J$-holomorphic disc $f_\infty:\D \to M$.
\end{lemma}

\proof Fix 3 distinct points $\zeta_j$, $ j=1,2,3$ on the unit circle and 3 distinct leaves $L_j$ of the characteristic foliation. Consider a sequence $\phi_k$ of Mobius transformations such that $f_k \circ \phi_k(\zeta_j) = L_j$ for $j=1,2,3$ and $k= 1,2,...$.  We claim that there are no disc bubbles for this sequence. By contradiction, a boundary disc bubble $g$  is attached to $S^2$ and everywhere transverse
to $\chi$ by part $(ii)$ of Proposition~\ref{translation}. Furthermore, by Gromov's Compactness Theorem (see \cite{mc-sa, ye}) the
boundary $g\vert_{\partial \D}$ of $g$ represents the trivial homological class
on $S^2 \backslash \{ p, q
\}$. 
Every $f_t(\partial \D)$ is a closed curve transverse to
$\chi$ and bounding a disc on $S^2$. 
Hence the curve $g(\partial \D)$ bounds on $S^2$ a totally
real disc $D$ (or a finite number of discs) and the line field $\chi$ does not vanish on its boundary. Then
it vanishes at an interior point of $D$ according to \cite{du-no-fo}, P.123, Cor. 14.5.2 : a contradiction. Q.E.D.

\vskip 0,1cm
{\it Step 4. Gluing Bishop's families.} The sequence $(f^{t_k})$
converges to a Bishop disc $f^{T}$. We point out that by the adjunction inequality (see \cite{mi-wh, mc}) the limit disc $f^T$ is embedded. If $f^T(\partial \D)$ does not intersect a sufficiently small neighborhood of $q$ on $S^2$ (to be precised below),
it is contained in the totally real part of $S^2$. Moreover according to Lemma~\ref{maslov-lemma} and to the fact that the winding number is a homotopy invariant, the winding number of $f^T$ is equal to 0.
It follows from Proposition~\ref{discdeformation-prop} that the disc $f^{T}$
generates a one-parameter family of nearby Bishop discs which contradicts the
maximality of $T$. Hence $f^T(\overline \D)$ intersects a sufficiently small neighborhood of $q$, meaning that it
intersects a boundary of some Bishop disc from the family $(g^t)$. By the
uniqueness Proposition~\ref{ellipticuniqueness-prop} the image of $f^T$ coincides with the image of one of the discs of this family. We recall that the Bishop discs as maps are defined up to three real parameters. Hence after a suitable reparametrization the two families of Bishop's discs are glued smoothly into a global one parameter family of discs with distinct images. This proves Theorem~\ref{elliptic-thm}. Q.E.D.

\section{Concluding remarks}

In this section we address several remarks concerning the results discussed in the present work.

{\bf Remark 1.}  For simplicity of notations we suppose everywhere that manifolds are $C^\infty$-smooth. However Theorem \ref{elliptic-thm} remains true under weaker assumptions. Indeed, it suffices to suppose that $J$ is of class $C^2$ and that the boundary of $M$ is of class $C^3$ for Theorem \ref{exhaust-theo} and other results, like the generation of discs near elliptic points or the transversality statements, to hold. The positivity of intersection in \cite{mi-wh} requires only $C^2$-smoothness of $J$.

{\bf Remark 2.} In Theorem \ref{elliptic-thm} it suffices to suppose that $\partial M$ contains no germs of non-constant $J$-holomorphic discs in a neighborhood of the sphere $S^2$. The proof does not require any change. We point out that in the smooth category this assumption is weaker than finite type assumptions. However, they are essentially equivalent in the real analytic case. See a more detailed discussion in \cite{ba-ma}. Using Theorem 
\ref{exhaust-theo} it is easy to prove (see \cite{di-su}) that if $\partial M$
is Levi convex and if there exists a strictly pseudoconvex function in a
neighborhood of $\partial M$, then $M$ can be exhausted by a sequence of
domains $M_j$ with strictly Levi-convex boundaries. Under this global
assumption the condition that $\partial M$ contains no non-constant
holomorphic discs can be dropped. Namely, we can consider a sequence $S^2_j$
of spheres in $\partial M_j$, their Levi-flat fillings, and pass to the
limit. This gives a  Levi-flat filling of $S^2$. However, if $\partial M$
contains non-constant holomorphic discs, the behaviour of boundaries of the
limit discs can be complicated as show examples from \cite{hi}.

{\bf Remark 3.} The assumption that $M$ contains no $J$-holomorphic spheres can be weakened. A simple topological argument of \cite{ye} shows that if a spherical bubble arises, then it has a negative self-intersection index. Thus it suffices to require that $M$ contains no sphere of this class.
One can also use a more subtil description of cusp-curves arising as bubbles (see \cite{ye}) and require that $M$ contains no such curves. This result is new for manifolds with Levi convex boundaries. However, this condition is difficult to verify. Finally, our method of proof allows  to obtain similarly  to \cite{ye} a Levi-flat filling with singularities if $M$ contains holomorphic spheres and $\partial M$ is Levi-convex. Thus Theorems 4,5,6 and 7 from \cite{ye} obtained there for strictly Levi convex manifolds remain true in the case where the boundary of $M$ is Levi convex. We point out that these results are new in the Levi-convex case. We do not state them here since they require a rather long description of generic assumption on an almost complex structure $J$ (rational regularity in the terminology of \cite{ye})  irrelevant to our forthcoming study of spheres with hyperbolic points. We leave the details to an interested reader.

\end{document}